\documentclass[a4paper,fleqn]{cas-sc}

\usepackage{filecontents}

\usepackage[numbers]{natbib}
\usepackage{multirow} 
\usepackage{amsmath}
\usepackage{amssymb}
\usepackage{algpseudocode}
\usepackage{mathtools}
\usepackage{algorithm} 
\usepackage{lineno}

\newtheorem{theorem}{Theorem}[section]
\newtheorem{lemma}[theorem]{Lemma}

\newtheorem{proposition}[theorem]{Proposition}
 \newproof{pf}{Proof}

\newdefinition{definition}{Definition}
\newdefinition{example}{Example}
\newdefinition{remark}{Remark}

\newcommand{\R}{\mathbb{R}}
\newcommand{\N}{\mathbb{N}}

\newcommand{\be}{\begin{equation}}
\newcommand{\ee}{\end{equation}}

\def\tsc#1{\csdef{#1}{\textsc{\lowercase{#1}}\xspace}}
\tsc{WGM}
\tsc{QE}
\tsc{EP}
\tsc{PMS}
\tsc{BEC}
\tsc{DE}

\begin{document}

\let\WriteBookmarks\relax
\def\floatpagepagefraction{1}
\def\textpagefraction{.001}
\shorttitle{Finite determination of accessibility and geometric structure of singular points for nonlinear systems}
\shortauthors{M. A. Sarafrazi et~al.}

\title [mode = title]{Finite determination of accessibility and geometric structure of singular points for nonlinear systems} 

%

\author[1]{Mohammad Amin Sarafrazi}[type=editor, 
        auid=000,bioid=1,
        orcid= 0000-0002-7265-0841]
\cormark[1]

\address[1]{ Rezvan complex,
Motahari Sq., Motahari Blvd.,  Shiraz 71868-98544, Iran}

\ead{sarafrazi@ut.ac.ir}



\author[2]{{\"U}lle  Kotta}[type=editor, 
        auid=000,bioid=1,
        orcid=0000-0002-7697-769X]

\address[2]{Department of Software Science, Tallinn University of Technology, Tallinn 12618, Estonia}

\ead{kotta@cc.ioc.ee}

\author[3]{Zbigniew Bartosiewicz}[
        auid=000,bioid=1,
        orcid=0000-0002-7250-4350]

\address[3]{Bialystok University of Technology, Faculty of Computer Science, Wiejska 45A, 15-351 Białystok, Poland}
\ead{z.bartosiewicz@pb.edu.pl}


\cortext[cor1]{Corresponding author}


\begin{abstract}
Exploiting tools from algebraic geometry, the problem of finiteness of determination of accessibility/strong accessibility is investigated for  polynomial systems and also for analytic systems that are immersible into polynomial systems. The results are constructive, and  algorithms are given to find the maximum depth of Lie brackets necessary for deciding accessibility/strong accessibility of the system at any point,  called here accessibility/strong accessibility index of the system, and is known as the degree of non-holonomy in the literature. Alternatively, upper bounds  on the accessibility/strong accessibility index are obtained, which can be computed easier.  In each approach, the entire set of accessibility/strong accessibility singular points are obtained.
Several examples demonstrate the applicability of the results using computer  algebra tools.

\end{abstract}




\maketitle

\section{Introduction}
Accessibility and strong accessibility are important notions in control theory, and necessary for most  control strategies. They are closely related to controllability, and in driftless systems or linear systems  become  equivalent  to controllability. Accessibility from a point of state space means the possibility of accessing an open set in the state space from  that point, using all possible inputs.

Similarly to controllability rank test for linear systems, there exists accessibility rank test for nonlinear systems. 
For analytic systems it is known that the control system is accessible from a point $x_0$ if and only if the dimension of the accessibility distribution $ C $ at this point is equal to the state dimension \cite{sussmann}.
 But, unlike the controllability of linear systems, in  nonlinear systems different points of the state space may have different accessibility properties, and as the accessibility distribution consists of infinite number of vector fields (Lie brackets of any depth), one does not know  to what extent the successive Lie brackets  need to be computed to make sure that the system is accessible/strongly accessible from a given point $x_0$, or $x_0$ is a singular point, i.e. the system is not accessible/strongly accessible from $x_0$  \cite{kawski}.
 
Similar problem exists in the context of controllability of non-holonomic systems, where the minimum depth of Lie brackets in the associated Lie algebra that determines controllability of the non-holonomic system at a point of the state space, called \emph{degree of non-holonomy} of the system at that point, is not known a-priori (see \cite{gabriel}, \cite{gabrielov}). It has been shown that for polynomial systems the maximum degree of non-holonomy over the entire state space is finite, however, to the best of our knowledge, no results are reported on the computation of the exact value of this maximal integer. In  \cite{risler} the author obtained an upper bound on the degree of non-holonomy for planar systems, and this result was extended to the case of polynomial systems with arbitrary dimension in  \cite{gabrielov}. The obtained upper bounds were improved and extended to the Noetherian analytic rings (i.e. a ring of analytic functions that is generated finitely) in \cite{gabrielov2} and \cite{gabriel}.
See also \cite{binyamini} for recent improvements. Unfortunately, these upper bounds grow drastically with increase in dimension of the system and degree of polynomials, and as a result, they are far from being applicable. 
 Singularity of distributions is also important in non-holonomic robots, where a configuration is called singular if the number of infinitesimal first order movements to reach nearby configurations increases  compared to other neighbor configurations \cite{muller, tchon}. 
 
 Knowing the exact location of the set of accessibility singular points, denoted by $S_\infty$ in this paper, is equally important, since $S_\infty$ is an invariant set, and therefore it should be avoided for initialization of the system, or trapping the state within  $S_\infty$. Also $S_\infty$ may obstruct global controllability, if the set of regular points is disconnected by the set $S_\infty$, just as singular points in flatness property can obstruct definition of global flat outputs and global motion planning\citep{kaminski}.
  
The approach presented in this paper, suggests that for  deciding  accessibility of a point in a finite number of steps, instead of examining accessibility property pointwise, one should look at  the big picture of the entire set of singular points and the invariance relations between them.
 The main idea  is the following: If the system is non-accessible from $x_0$, then any trajectory starting from $x_0$ must evolve on the invariant set of non-accessible points. It is shown that in polynomial systems, and analytic systems that are immersible into polynomial systems, the set of singular points of accessibility  are algebraic sets. Then   the invariance of  algebraic sets are characterized in terms of invariance of their corresponding ideals under Lie derivations defined by the system dynamics.
The polynomial structure of the vector fields that describe the system is responsible for stabilization of the constructed sequences of sets that at the limit gives us the invariant set. Note that similar sequences for analytic or meromorphic vector fields do not have to stabilize. Analogous results have been obtained for the case of strong accessibility. \\
Our result provides finite and applicable accessibility tests in several different ways.  A group of approaches presented in this paper gives either the exact or the upper bound on the depth of Lie brackets of vector fields in the (strong) accessibility distribution that one needs to compute for the usual accessibility/strong accessibility rank test.  Simultaneously, in  all the approaches the entire set of singular points is obtained as the algebraic set of a limiting ideal of an ascending chain of ideals that stabilizes, and its stabilization can be detected constructively by a  differential algebraic test. 
In each of the  proposed methods, only one chain of ideals  is sufficient to determine the singular points of the entire state space. \\
 We restrict the main results  on polynomial systems, because of Noetherian property of the ring of polynomials, and especially because it is easier to manipulate ideals of the ring of polynomials using  computer algebra tools. We then extend the results to the case of analytical systems that can be immersed into polynomial systems, which contains a  very wide class of systems. 
  
  Finally, note that similar results have been obtained in \cite{amin} based on the same idea, for rational discrete-time systems, and also for analytic systems restricted to a compact semianalytic set,
where it has been shown that the set of singular points of accessibility is the limiting algebraic set of a specific descending chain of algebraic sets $S_k$, with $S_k$ being the set of states from which the system is not accessible in $k$ steps. Similarly, it has been shown in \cite{amin} that a certain integer $r^*$, named accessibility index of the system,  can be found such that for any point of state space, the discrete-time system is  accessible if and only if it is  accessible for input sequences of length $r^*$, and hence renders the infinite  aceessibility test to a finite test. However, in the discrete-time case it is possible to compute explicitly the solution of state evolution at any time instance $k$, which gives a simple characterization of set of accessibility singular points, and results in a strictly descending chain of algebraic sets $S_k$.
  A similar approach for the continuous-time case would inevitably lead to complications as we are not able to compute solutions and the appropriate sets, now parametrized by the continuous time $t$, and the analogous chain of algebraic sets $S_k$ (see Definition \ref{def5}) may not be strictly descending. Therefore a different characterization of the set of accessibility singular points is needed.

The paper is organized as follows. Preliminaries and definitions are given in Section \ref{sec2}.
  Generic accessibility  criterion, and
relation between generic accessibility and pointwise accessibility
are obtained in Section \ref{sec3}. The main results of the paper for polynomial systems are given in Section \ref{sec4}. Section \ref{sec5} contains the extension of the results of the polynomial case to the case of non-polynomial  systems. Appendix presents an introduction on ideals and algebraic sets.
 
 \section{Preliminaries} \label{sec2}

We denote by $\mathcal{A}$ the set of analytic functions of $x$ on $\R^n$, and by $\R [x]= \R [x_1 , \cdots , x_n] $ the commutative ring of all polynomials in variables $x_1,...,x_n$  with coefficients in $\R$.
By $\R[x]^n$ we denote the set of all vector fields on $\R^n$ with components in $\R[x]$.  \\
Let $V^{\omega}$ be the Lie algebra of analytic vector fields on $\R^n$. For  $f,g\in V^{\omega}$ and  $\varphi \in \mathcal{A}$, let $L_f \varphi \coloneqq (\frac{\partial \varphi}{\partial x})f$  be the Lie derivative of the function $\varphi$ along the vector field $f$, and $[f,g]\coloneqq \frac{\partial g}{\partial x}f - \frac{\partial f}{\partial x}g$
 be the Lie bracket of $f,g$. Also we use the notation $ad_f(g):=[f,g]$.  For $f \in \R[x]^n$ and $\varphi \in \R[x]$, we have $L_f \varphi \in \R[x]$ and for $f,g \in \R[x]^n$ we have $[f,g] \in \R[x]^n$. Obviously $\R[x]^n$ is a Noetherian module over  $\R[x]$, and therefore any submodule of $\R[x]^n$ is finitely generated. Moreover, any ascending chain of submodules of $\R[x]^n$ eventually stabilizes \cite{zariski}. \\
Recall that for   $f,g \in V^{\omega}$, scalar   $r\in \R, $   and functions $p_1, p_2 \in \mathcal{A}$, the following properties hold:
 \begin{eqnarray}
&& \hskip -1cm[f,g]=-[g,f],  \label{1a} \\ 
 &&\hskip -1cm [r f_1 + f_2, g ]= r [f_1 ,g] +[f_2 ,g]   \label{1b} \\ 
 &&\hskip -1cm [p_1 f,p_2g]=p_1 p_2 [f,g]+(L_f p_2)p_1 g -(L_g p_1 )p_2 f    \label{1c}
 \end{eqnarray} 
A (analytic) distribution $\cal{D}$ assigns to each point $x\in \R^n$ a linear subspace of the tangent space $\R^n$. We say that a distribution $\cal{D}$ is generated by a set of vector fields $\{ f_1,\dots,f_k\}$ if $\mathcal{D}(x)=\textup{span}_{\R}\{ f_1(x),\dots,f_k(x)\}$ at every $x$, and in this case we identify the distribution $\cal{D}$ by its generators.
A distribution $\cal{D}$ is said to be invariant under a vector field $Y$ if $[X,Y] \in \cal{D}$ whenever $X \in \cal{D}$.
Consider the nonlinear  system described by the equation of the form
\begin{equation} \label{sys}  \Sigma : \dot{x}(t)=f(x(t))+ \sum_{i=1}^m u_i(t)g_i(x(t))
\end{equation}
where $t \in \R,~ x(t)\in D \subset \R^n , u(t)=(u_1(t),...,u_m(t))\in U \subset\R^m$, and $f,g_1,\dots,g_m$ are analytic vector fields. The sets $D$ and $U$ are assumed to be open. 
We denote by $\mathcal{U}$ the set of all measurable locally essentially bounded maps $u:[0,\infty) \rightarrow U$,
and by  $\Phi^{\Sigma ,u}_{x_0} (t)$  the trajectory of the system $\Sigma$ at time instant $t$, starting from the initial state $x_0 \in D$ and driven by the input function $u \in \mathcal{U}$.
The set of reachable points from $x_0$ at time (exactly) $t>0$ is denoted by $\mathcal{R}_{\Sigma}(x_0,t)\coloneqq \{ \Phi^{\Sigma ,u}_{x_0} (t) ~|~ u \in \mathcal{U}\}$, 
and the set of points reachable from $x_0$  is denoted by   $\mathcal{R}_{\Sigma}(x_0)\coloneqq \bigcup_{t\geq 0} \mathcal{R}_\Sigma (x_0,t)$.
\begin{definition}
\cite{sussmann} The system $\Sigma$ is said to be \emph{accessible from} $x_0$ if  $\textup{int}(\mathcal{R}_{\Sigma}(x_0))\neq \emptyset$. The system $\Sigma$ is said to be \emph{strongly accessible from} $x_0$  if $\textup{int}(\mathcal{R}_{\Sigma} (x_0,t))\neq \emptyset$ for every $t > 0$. A point $x^*$ is called a singular point of accessibility (strong accessibility) for the system $\Sigma$ if the system $\Sigma$ is not accessible (strongly accessible) from $x^*$. 
 \end{definition}
For analytic systems, the above definition of strong accessibility is equivalent to  $\textup{int}(\mathcal{R}_{\Sigma} (x_0,t))\neq \emptyset$ for every $0 < t \leq T$ for some $T>0$ (see \cite{sussmann}).
 \begin{definition}
\cite{nijmeijer} Consider the nonlinear system (\ref{sys}). The \emph{accessibility algebra} $\mathcal{C}$  is the smallest subalgebra of $V^{\omega}$ that contains $\{f,g_1,...,g_m \}$, and the \emph{accessibility distribution} $C$ is the distribution generated by the accessibility algebra $\mathcal{C}$.  The \emph{strong accessibility algebra}  $\mathcal{C}_0$   is the smallest  subalgebra of $V^{\omega}$ that contains $\{g_1,...,g_m \}$ and is invariant under $ad_f$, and the \emph{strong accessibility distribution} $C_0$ is the distribution generated by the strong accessibility algebra $\mathcal{C}_0$.   
 \end{definition}
Every element of $\mathcal{C}_0$ is a linear combination of repeated Lie brackets of the form $[X_k,[X_{k-1},[\cdots,[X_1,g_j]\cdots ]]], ~~1\leq j \leq m, k\in \N$ and $X_i \in \{f,g_1,\dots,g_m\}$, and $\mathcal{C}=\mathcal{C}_0\cup f$. Both $C$ and $C_0$ are involutive distributions.
 
The accessibility and strong accessibility  can be determined by the so-called accessibility rank condition:
\begin{theorem} \label{theorem1}
\cite{sussmann} The system (\ref{sys}) is accessible (respectively strongly accessible) from a point $x_0 \in D$ if and only if $\textup{dim}\! ~C(x_0)$ $=$ $n$ (respectively $\textup{dim}\!~ C_0(x_0)=n$). 
\end{theorem}

Thanks to involutivity,   $C$    has maximal integral submanifold property \cite{sussmann2}, which means that through every point $x_0\in D$ passes a (unique) maximal integral submanifold  $I(C,x_0)$, such that for  every $x\in I(C,x_0)$, the  tangent space of  $ I(C,x_0)$ at $x$ is equal to $C(x)$. Each  $I(C,x_0)$ is a forward-invariant set for the system.
The set  $\mathcal{R}_{\Sigma}(x_0)$ is contained in  $ I(C,x_0)$ and has nonempty interior in it.

Now we define a filtration of  accessibility (respectively strong accessibility)  distributions of order $k$, as well as a descending chain of algebraic sets, corresponding to singular points of each distribution. Our main result characterizes the limiting algebraic sets in terms of invariance with respect to the system vector fields.
\begin{definition} \label{def5}
For $k\geq 0$ we denote by $\mathcal{C}^k$ (respectively $\mathcal{C}_0^k$) the  smallest subset  of $\mathcal{C}$ (respectively $\mathcal{C}_0$) that contains all Lie brackets of depth at most $k$ from the accessibility algebra $\mathcal{C}$ (respectively  the strong accessibility algebra $\mathcal{C}_0$), and correspondingly define \emph{accessibility distribution of order $k$}, denoted by $C^k$ (respectively \emph{strong accessibility distribution of order $k$}, denoted by $C_0^k$) as the distribution generated by it.
 We denote by $S_k$ (respectively $S_k^*$) the set of all points $x \in D$  such that $\textup{dim}\!~C^k(x)<n$ (respectively $\textup{dim}\!~C_0^k(x)<n$), and by $S_{\infty}$ (respectively $S_{\infty}^*$) the set of points $x\in D$ such that $\textup{dim}\!~C(x)<n$ (respectively $\textup{dim}\!~C_0(x)$ $<$ $n$). By Theorem \ref{theorem1} the set $S_{\infty}$ (respectively $S_{\infty}^*$) is the set of singular points of accessibility (respectively strong accessibility).
\end{definition}



\section{Generic  versus pointwise properties} \label{sec3}
Recall that a property is said to hold generically if it holds \emph{almost} everywhere, i.e. except on a set of measure zero. For an  analytic distribution $\mathcal{D}$, due to  analyticity, the generic dimension  is  the maximum dimension it can have at any $x\in D$. 
We show that generic accessibility (respectively generic strong accessibility) of the system is necessary for accessibility (respectively strong accessibility) from every individual point, and provide a finite test for checking this property. Therefore we single out those systems that are not generically accessible (respectively generically strongly accessible).

To avoid confusion, by $\textup{dim}\!~\mathcal{D}(x)$ we mean the dimension of $\mathcal{D}$ evaluated at the point $x$, while we use $\textup{dim}\!~\mathcal{D}$ to denote the generic dimension of $\mathcal{D}$ over all $x\in D$.

\begin{theorem} \label{theorem33}
The analytic system (\ref{sys}) is generically accessible (respectively generically strongly accessible) if and only if    $\textup{dim}\!~C^{n-1}=n$ (respectively $\textup{dim}\!~C_0^{n-1}=n$). Moreover, if   $\textup{dim}\!~C^{n-1}\neq n$ (respectively $\textup{dim}\!~C_0^{n-1}\neq n$), then the system is non-accessible (respectively strongly non-accessible) from every $x\in D$. 
\end{theorem}
\begin{pf}
Consider the chain of distributions $C^0\subset C^1\subset \cdots$. Since $V^{\omega}$ is an $n$-dimensional vector space, therefore for some $k^*\leq n-1$ we have $\textup{dim}\!~C^{k^*}=\textup{dim}\!~C^{k^*+1}$.
Assume that  $\{ h_1,\dots, h_q \}$ are vector fields from $\mathcal{C}^{k^*}$ that generate $C^{k^*}$.
 By construction, $C^{k^*+1}$ is generated by vector fields $\{ h_1,\dots, h_q \}$ together with vector fields of the form $ad_X{h_i}$  for all $1\leq i \leq q$ and all $X\in \{f,g_1,\dots, g_m\}$. 
Therefore the assumption  $\textup{dim}\!~C^{k^*}=\textup{dim}\!~C^{k^*+1}$ means that   $ad_X{h_i}$ $\in$ $\textup{span}_{\mathcal{A}}\{h_1,\dots, h_q\}$, for any $1\leq i \leq q$ and any $X$ $\in$ $\{f,g_1,\dots, g_m\}$.  By (\ref{1b}) and (\ref{1c}) and a simple induction it follows that $ ad_{X_1}ad_{X_2}\dots ad_{X_j}(h_i)$ $\in$ $\textup{span}_{\mathcal{A}}\{h_1,\dots, h_q\}$ for  any  $X$ $\in$ $\{f,g_1,\dots, g_m\}$,  $j \in \N$, and $1\leq i \leq q$.
Hence, for $k\geq k^*$, we have $\textup{dim}\!~C^{k}=\textup{dim}\!~C^{k^*}$, which gives  $\textup{dim}~\!C^{n-1}=\textup{dim}\!~C$, and therefore the system is generically accessible if and only if $\textup{dim}\!~C^{n-1}=n$. 
Also, since the generic dimension of an analytic distribution is the maximum dimension it can have, therefore if $\textup{dim}\!~C^{n-1}<n$, then at every individual point $x\in D$ we have $\textup{dim}\!~C(x)<n$.\\
 By replacing $C^i$ in the above with $C_0^i$, the proof of strong accessibility case follows similarly.    \qed
\end{pf} 
\begin{definition}  \label{MNAD}
For a generically accessible  (respectively generically strongly accessible) system, the  integer $r^*$ (respectively $l^*$) is called the \emph{accessibility index} (respectively \emph{strong accessibility index}) of the system $\Sigma$ over $D$ if it is the maximum integer for which there exists at least one point
$x_0 \in D$ such that $\textup{dim}\!~C^{r^*-1}(x_0)<n$ and $\textup{dim}\!~C^{r^*}(x_0)=n$ (respectively $\textup{dim}$ $\!~C_0^{l^*-1}(x_0)$ $<$ $n$ and $\textup{dim}\!~C_0^{l^*}(x_0)=n$). If  there is no such finite integer $r^*$ (respectively $l^*$), we put $r^*=\infty$ (respectively $l^*=\infty$).
 \end{definition}


 \begin{remark} \label{remark1}
In the literature on differential geometry, the singular points of distribution have been studied in the context of singular foliation theory
 \cite{sussmann2,nagano}, 
where the integral manifolds of a distribution are seen as leaves of a foliation, and when the leaves have not the same dimension, we have a singular foliation. 
From  definition, the set $S_\infty$ is the union of singular leaves (leaves that are of lower dimension with respect to neighbor leaves), i.e. every singular leaf is contained in $S_\infty$, and also through every $x_0 \in S_\infty$  passes a singular leaf. 
But, some care should be paid in distinguishing between $S_\infty$ and singular integral manifolds. The set $S_\infty$ is an algebraic set, and not necessarily a manifold.  Even in the case when $S_\infty$ is a manifold,  the geometric dimension of $S_\infty$ may be greater than the dimension of  contained integral manifolds,  as the following example shows.
\end{remark}
\begin{example} \label{ex51}
Consider the system defined by $\Sigma: \dot{x}=f+ug$ with $f:=[0~~x_2^2+x_3^2-1~~0]^T$ and $g:=[x_2~~x_2x_3~~-x_2^2]^T$. It is easy to verify that $\Sigma$ is generically accessible, with the set of accessibility singular points $S_\infty=\{x\in \R^3 | ( x_2^2+x_3^2-1)=0\}$ (using Theorem \ref{theorem313} bellow). Although the geometric dimension of $S_\infty$ is 2, for every point $x\in S_\infty$ the accessibility distribution is one dimensional, which means that the integral manifold of such points is a line,  described by the non-polynomial equations $x_2=\sin (x_1), ~ x_3=\cos (x_1)$. 
\end{example}

\section {Polynomial systems} \label{sec4}

\subsection{Singular points of accessibility} \label{3.1}

\begin{lemma} \label{lemma2}
For a polynomial system $\Sigma$, the accessibility index $r^*$ is finite, and the set of singular points of accessibility is an algebraic set.

\end{lemma}
\begin{pf}
Denote by $M_k$ a matrix whose columns are vector fields  of the set $\mathcal{C}^k$, defined in Definition \ref{def5}. Assume that $M_k$ has $l$ minors of dimension $n\times n$, denoted by $m_{k,1},...,m_{k,l}$. By Definition \ref{def5}, the set $S_k$ is the set 
\begin{equation} \label{2a}
S_k=\{x\in D~ |~ m_{k,i}(x)=0~ \textup{for all}~1\leq i \leq l \} 
\end{equation}
Since all vector fields are assumed to be polynomial vector fields, (\ref{2a}) determines an algebraic set [3], corresponding to the ideal $I_{M_k}\coloneqq \langle m_{k,1},...,m_{k,l} \rangle$, or say, $S_k =\mathcal{V}(I_{M_k})$ (see Appendix for definition of the operator $\mathcal{V}$). Hence every set $S_k$ is an algebraic set.
Now, by construction, for any $k$, the matrix $M_k$ is a submatrix of $M_{k+1}$. Thus all minors of $M_k$ are minors of $M_{k+1}$ too, and therefore $I_{M_k}\subseteq I_{M_{k+1}}$, from which, using the Proposition 1 in Appendix we have $S_k \supseteq S_{k+1}$. Hence we have the following descending chain of algebraic sets
\begin{equation} \label{3}
S_1 \supseteq S_2 \supseteq \cdots
\end{equation}
and because the ring $\R [x_1 ,...,x_n]$ is a Noetherian ring, from Hilbert Basis Theorem [3] it follows that the chain (\ref{3}) eventually stabilizes, i.e.,  there exists some integer $r^*$ such that 
\begin{equation} \label{4}
S_1 \supseteq S_2 \supseteq \cdots  \supseteq S_{r^*} =S_{r^*+1}= \cdots = S_{\infty}
\end{equation}
The smallest integer $r^*$ satisfying  (\ref{4}), by definitions \ref{def5}  and \ref{MNAD},  is the accessibility index of the system $\Sigma$. 
Also, the set  $S_{\infty}=S_{r^*}$ is an algebraic set. \qed
\end{pf} 

Lemma \ref{lemma2} states that the polynomial system $\Sigma$ has finite accessibility index, but,  unfortunately, it is based on the Hilbert Basis Theorem, which is not a constructive theorem, and doesn't warrant the inclusion relations in (\ref{4}) to be exclusive. In other words, it is not clear when the chain of algebraic sets stabilizes forever.  In the following, our main result addresses this problem. 
%

 \begin{theorem}   \label{theorem3}
Consider  a generically accessible system of the form (\ref{sys}) and the set $S_\infty$ as defined in Definition \ref{def5}. Then $S_\infty$ is the maximal zero-measure forward-invariant  set of the system.
\end{theorem}
\begin{pf}
 By Lemma \ref{lemma2}, the set $S_\infty$ for a generically accessible system is a closed zero-measure set. 
Every point $x_0 \in S_\infty$ is contained in a maximal integral manifold $I(C,x_0)$ which is a forward invariant set for the system, and every other point $x\in I(C,x_0)$  belongs to $S_\infty$ Therefore $S_\infty$ is a zero-measure forward-invariant set for the system.
Let  $B$ be any  zero-measure forward-invariant set of the system. For every $x_0 \in B$, the set of reachable points from $x_0$ is contained in the zero-measure set $B$, and hence it has empty interior. This means that all points of $B$ belong to $S_\infty$.
\end{pf}
\begin{definition}  \label{def7}
For an ideal $I$ of $\R[x]$ and a polynomial vector field $f \in  \R[x]^n$,
we say that $I$ is invariant under the operator $L_f$ if for every $p \in I$ we have $L_f (p) \in I$.
\end{definition}
\begin{proposition}   \label{theorem14}
 For an ideal $I=\langle z_1,...,z_k \rangle$ of  $\R[x] $ and a vector field $X\in \R^n[x]$, the ideal   $I$ is invariant under   $L_X$ if and only if $L_X(z_i)\in I$ for $1\leq i \leq k$.
\end{proposition}

In what follows, we use $\mathcal{I}(A)$ to denote the zero-ideal of a given set $A\subset \R^n$ (see Appendix for a formal definition of zero-ideal of a set).
\begin{lemma} \label{lemma38}
For a polynomial system of the form (\ref{sys}), an algebraic set $V$ is  forward-invariant  if and only if for every $X\in \{f,g_1,...,g_m\}$, the ideal $\mathcal{I}(V)$ is invariant under  $L_X$.
\end{lemma}
\begin{pf}
\emph{Sufficiency.} Assume $\mathcal{I}(V)$ to be 
invariant under $L_X$, for every $X\in \{f,g_1,...,g_m\}$.
 First we show that the set $V$ is forward-invariant under any constant input. From the assumption we get that for any $p\in \mathcal{I}(V)$, we have $L_Y(p)\in \mathcal{I}(V)$ for any $Y=f+\sum_{i=1}^m c_i g_i$, where $c_1,...,c_m$ are constants. Inductively, we have $L_{Y}^{(r)}(p)\in \mathcal{I}(V)$ for any $r\in \N$. Assume that $\mathcal{I}(V)$ is generated by the polynomials $\{z_1,...,z_l\}$. Now, because $z_1,...,z_l$ belong to the ideal $\mathcal{I}(V)$, we conclude that for any $z_i$, we have  $L_{Y}^{(r)}(z_i)\in \mathcal{I}(V)$ . Because at the time $t=0$ we have
\begin{equation}
\frac{d^r}{dt^r}z_i(x(t)) \Bigm\vert_{t=0} =L_{Y}^{(r)}(z_i)(x(0)),
\end{equation}
we conclude that for an initial state $x\in V$, and under the constant input $u=(c_1,...,c_m)$, the functions $z_1,...,z_l$ remain zero along the trajectory of the system. In other words, the trajectory lies completely in $V$. 
Inductively, we can conclude that for any piecewise constant input, the set $V$ is forward-invariant. 
This result can be extended from piecewise constant inputs to any input $u\in \mathcal{U}$, because the set of piecewise constant inputs is dense in $\mathcal{U}$, meaning that for any input  $u\in \mathcal{U}$, any $x_0\in D$, any time $t$ and any $\epsilon > 0$, one can find a piecewise constant $\bar{u}$ such that $\vert \Phi_{x_0}^{\Sigma,u} (t) -\Phi_{x_0}^{\Sigma,\bar{u}}(t) \vert < \epsilon $ (see Lemma 2.8.2 in [4]), so $V$ is forward-invariant under any $u\in  \mathcal{U}$.
\\
\emph{Necessity.} Assume that $\mathcal{I}(V)$ is generated by the polynomials $\{z_1,...,z_l\}$ and is not invariant under all $L_X$, for every $X\in \{f,g_1,...,g_m\}$. Then by Proposition \ref{theorem14}, for some $1\leq i \leq l$ and some $X\in \{f,g_1,...,g_m\}$ we have $ L_X(z_i) \notin \mathcal{I}(V)$. This means that there is some point $x_0 \in V$ such that $L_X(z_i)(x_0) \neq 0$. This implies that the trajectory of the state starting from $x_0$ in the direction of the vector field $X$ does not lie completely in $V$. In fact, for every trajectory $x(t)=\Phi_{x_0}^{\Sigma,u} (t)$ that lies completely in $V$, we must have $z_i(x(t))=0$ for all $t\geq 0$ and all $1\leq i \leq l$, and therefore
\begin{equation}
\frac{d}{dt}z_i(x(t))\Bigm\vert_{t=0} =L_X(z_i)(x_0)=0~~ \textup{for all}~~ 1 \leq i \leq l.
\end{equation}
So $V$ cannot be a forward-invariant set of the system. \qed
\end{pf}
\begin{theorem} \label{theorem39}
Consider a polynomial system $\Sigma$ of the form (\ref{sys}) and the sets $S_k$ as defined in Definition \ref{def5}. Then
\begin{enumerate}
\item[(a).] The ideal $\mathcal{I}(S_{\infty})$ is the minimal  real radical ideal that is invariant under $L_f,L_{g_1},$ $\dots,L_{g_m}$.
\item[(b).] The accessibility index of  $\Sigma$ is the smallest  $k$ for which $\mathcal{I}(S_k)$ is proper and invariant under $L_f,L_{g_1},$ $\dots,$ $L_{g_m}$. 

\end{enumerate}
\end{theorem}
\begin{pf}
\emph{(a).}
By Lemma \ref{lemma2},  $S_\infty$ is an algebraic set. Then   from Theorem \ref{theorem3} and Lemma \ref{lemma38} the claim is obvious.

\emph{(b).} Let $\mathcal{I}(S_k)$  be proper and invariant under $L_f,L_{g_1},$ $\dots,$ $L_{g_m}$. Based on Lemma \ref{lemma38}, the set $S_k$ must be an invariant set of the system, which from Theorem \ref{theorem3} gives $S_k \subset S_\infty$. On the other hand, from (\ref{4}) we have $ S_\infty  \subset S_k $ for any $k$. Therefore  $ S_\infty  = S_k $. From (\ref{4}) in the proof of Lemma \ref{lemma2}, the accessibility index of the system is the smallest integer $r^*$ such that $S_{r^*}=S_\infty$, and the claim is proved.  \qed
\end{pf}

\begin{algorithm}[H]
\caption{(Computing $S_\infty$ and accessibility index $r^*$, assuming the generic rank of $M_{k^*}$ is full)}
\label{<your label for references later in your document>}
\begin{algorithmic}[1]
\State \textbf{Initialization:} $k \gets k^*,  M_k \gets M_{k^*}$
\State Compute all $n \times n$ minors of the matrix $M_k$
\State Construct the ideal $I_{M_k}$ (as in the proof of Lemma \ref{lemma2})
\State Compute $\sqrt[\R]{I_{M_k}}$
\State Compute a  basis $\{b_1 ,...,b_r\}$ for $\sqrt[\R]{I_{M_k}}$
\If {for ever $b_i$ we have $L_X(b_i)\in \sqrt[\R]{I_{M_k}}$ for all $X \in \{f,g_1 , ..., g_m \}$,}
    \State stop and return $S_\infty  \gets \mathcal{V}(I_{M_k})$ and $ r^* \gets k$
\Else
    \State $k \gets k+1$ and go to step 2
\EndIf
\end{algorithmic}
\end{algorithm}

Based on Theorem \ref{theorem39}, Algorithm 1  can be used for finding the set of accessibility singular points $S_\infty$ and the accessibility index of the system.

%
%

\begin{example} \label{example}
Consider the nonlinear control system
   \begin{equation}  \label{ex1} \left  \{ \begin{array}{lll}     \dot{x}_1=u_1 x_2 \\ \dot{x}_2=u_2 x_1^2    \end{array}  \right.   \end{equation}   
   We use Algorithm 1  to obtain the singular points of accessibility and accessibility index. For this system, we have $  g_1=[ x_2 ~~ 0]^T$  and $g_2=[ 0 ~~ x_1^2]^T $, and 
\be \label{M0ex1}       M_0(x_1,x_2)=[ g_1 ~~ g_2]=\begin{bmatrix} x_2  & 0 \\ 0  & x_1^2  \end{bmatrix} \ee
 The generic rank of the matrix $M_0$ is  2. So we initialize the Algorithm 1 with $M_0$. The ideal $I_{M_0}$ is generated by the determinants of all $2 \times 2$ minors of$M_0$:
\begin{equation} \label{IM0ex1} I_{M_0} =\langle x_1^2x_2 \rangle~ \rightarrow ~ \sqrt[\R]{I_{M_0}} =\langle x_1x_2 \rangle. \end{equation}
 We chack the invariance of  $ \sqrt[\R]{I_{M_0}}$ under $L_{g_1}$ and $L_{g_2}$. We have $ L_{g_1}(x_1x_2)=x_2^2 \notin \langle x_1x_2 \rangle $,  
which shows that it is not invariant. So, according to Algorithm 1, we proceed to the next step and obtain $M_1$ and $I_{M_1}$, and perform the previous computations again. We have $[g_1,g_2]=[ -x_1^2 ~~  2x_1x_2 ]^T$, and
\[  M_1(x_1,x_2)=\begin{bmatrix} g_1 & g_2 &  [g_1,g_2]  \end{bmatrix}=\begin{bmatrix} x_2 & 0 & -x_1^2 \\ 0& x_1^2 & 2x_1x_2 \end{bmatrix} \] 
The ideal $I_{M_1}$ is generated by the determinants of all $2\times 2$ minors of $M_1$. So we have $ I_{M_1}=\langle x_1^2x_2, x_1x_2^2, x_1^4\rangle$, and therefore  $\sqrt[\R]{I_{M_1}}=\langle x_1 \rangle$.
We check the invariance of  $\sqrt[\R]{I_{M_1}}$ under $L_{g_1}$ and $L_{g_2}$. We have $ L_{g_1}(x_1) =x_2 \notin  \langle x_1 \rangle.$
So we proceed by obtaining $M_2$ and $I_{M_2}$, and performing the previous computations again.
\[  M_2(x_1,x_2)=\begin{bmatrix} M_1 & [g_1,[g_1,g_2]] & [g_2,[g_1,g_2]]\end{bmatrix}\]
\[=\begin{bmatrix} x_2 & 0 & -x_1^2 & -4x_1x_2 & 0 \\ 0 & x_1^2 & 2x_1x_2 &  2x_2^2 & 4x_1^3 \end{bmatrix} \] 
 $I_{M_2}$ is generated by the determinants of all $2\times 2$ minors of $M_2$. So  $ I_{M_2} =\langle x_1^2x_2,x_1x_2^2,x_2^3,x_1^4 \rangle$, and $\sqrt[\R]{I_{M_2}} =\langle x_1,x_2 \rangle$. 
 We check whether the ideal is invariant under $L_{g_1}$ and $L_{g_2}$:
\begin{eqnarray} \nonumber
 &L_{g_1}(x_1)
 =x_2 \in \langle x_1,x_2 \rangle, ~  L_{g_2}(x_1)
=0 \in \langle x_1,x_2 \rangle, \\    \nonumber
& L_{g_1}(x_2)
=0 \in \langle x_1,x_2 \rangle ,~L_{g_2}(x_2)
=x_1^2 \in \langle x_1,x_2 \rangle. 
\end{eqnarray}
This shows that $\sqrt[\R]{I_{M_2}}$ is invariant under the vector fields of the system and therefore, according to Algorithm 1, the set $\mathcal{V}(I_{M_2})=S_\infty=(0,0)$ is the set of accessibility singular points,  and the accessibility index of the system is 2, which means that computation of Lie brackets of depth up to 2 determines accessibility for every point. For comparison, the results of \cite{risler} suggests that for a polynomial system of order 2 and degree of polynomials no more than $d$, the Lie brackets of depth up to $6d^2 -2d+2$ may be needed, which for this example means all Lie brackets of depth up to 22.
\end{example}


Computation of real radical for general ideals is a challenging task, and this motivates us to propose alternative approaches for obtaining the set of singular points, as well as upper bounds on the accessibility index, which can be computed easier.

\begin{theorem} \label{theorem313}
For a  polynomial system of the form (\ref{sys}), assume that $I_{M_q}$ is proper  for some $q<n$. Let   $\bar{I}_{M_{q}}$ be the smallest ideal that contains $I_{M_{q}}$ and is invariant under  $L_f,$ $L_{g_1},$ $\dots,L_{g_m}$. Then $\mathcal{V} (\bar{I}_{M_{q}})=S_\infty$.
\end{theorem}
\begin{pf}
Note that  Theorem \ref{theorem33} assures that for a generically accessible system,    $I_{M_{q}}$ is a proper ideal for some $q<n$.
From the proof of Lemma \ref{lemma2}, $S_q=\mathcal{V}(I_{M_q})$, which, using Proposition \ref{prop2} in Appendix, and part (ii) of Proposition \ref{prop1} in Appendix gives $ I_{M_{q}} \subset \mathcal{I}(S_{q})$. On the other hand, from (\ref{4}) and part (vi) of Proposition \ref{prop3} in Appendix, we have 
$\mathcal{I}(S_{q}) \subset \mathcal{I}(S_{\infty})$. Theses last two relations give
\be  \label{harchi}
\mathcal{I}_{M_q} \subset \mathcal{I}(S_{\infty}).
 \end{equation} 
Since by part (a) of Theorem \ref{theorem39} the ideal $\mathcal{I} (S_\infty)$ is invariant under $L_f,L_{g_1},\dots,L_{g_m}$, therefore from (\ref{harchi}) and the definition of $\bar{I}_{M_{q}}$ we get that $\bar{I}_{M_{q}}\subset \mathcal{I} (S_\infty)$. Therefore
\begin{equation} \label{harchi2}
S_\infty \subset \mathcal{V} (\bar{I}_{M_{q}}).
\end{equation}
On the other hand, from the definition of $\bar{I}_{M_{q}}$ and Lemma \ref{lemma38} we  conclude that $\mathcal{V} (\bar{I}_{M_{q}})$ is a forward-invariant set of the system, and since $\bar{I}_{M_{q}}$ is a proper ideal, the set $\mathcal{V} (\bar{I}_{M_{q}})$ is a zero-measure  set. Therefore from Theorem \ref{theorem3} we have $\mathcal{V} (\bar{I}_{M_{q}}) \subset S_\infty$, which together with (\ref{harchi2}) gives $S_\infty = \mathcal{V} (\bar{I}_{M_{q}})$.  \qed
\end{pf}

For a given ideal $I=\langle p_1,...,p_r \rangle$, to compute $\bar{I}$ as in Theorem \ref{theorem313}, it suffices to apply the operators  $L_f,$ $L_{g_1},$ $\dots,L_{g_m}$ to the generators of $I$, and then constitute  a new ideal generated by  $ \{p_1,...,p_r\}$ and   $\{L_f(p_j),$ $L_{g_1}(p_j),$ $\dots,$ $L_{g_m}(p_j)$, $1\leq j \leq r \}$, and continue this procedure inductively,  until the new ideal becomes equal to the previous one. The stabilization of this procedure is guaranteed by the Hilbert Basis Theorem \cite{cox}. 

 Theorem \ref{theorem313} leads to another algorithm (Algorithm 2 ) for obtaining the entire set $S_\infty$, without the need for computing real radical of ideals. 
 \begin{algorithm}[H] 
\caption{(Computing $S_\infty$, assuming the generic rank of  $C^q$ is complete)}
\label{<your label for references later in your document>}
\begin{algorithmic}[1]
\State Compute $I_{M_q}$ (as in the proof of Lemma \ref{lemma2})
\State  $\bar{I} \gets I_{M_q}$
\State $J \gets \emptyset$
\State Assuming $\bar{I}=\langle z_1,\dots,z_k \rangle$ :
\For{i:1,\dots,k}
  \For{every $X_j\in \{f,g_1,\dots,g_m\}$}
 \State $ J \gets J \cup \langle L_{X_j}(z_i) \rangle$
  \EndFor
\EndFor
\If {$J\subset \bar{I}$}
    \State Stop and return $S_\infty \gets  \mathcal{V}(\bar{I})$ 
\Else
    \State $\bar{I} \gets \bar{I}\cup J$ and go to step 3
\EndIf
\end{algorithmic}
\end{algorithm}

\begin{example}
Let us consider the system (\ref{ex1}) of Example \ref{example}  again and obtain $S_\infty$ for this system using Algorithm 2. It can be seen from  (\ref{M0ex1}) that the generic rank of  $C^0$ is 2. So we use  Algorithm 2, with  $I_{M_0} =\langle x_1^2x_2 \rangle$ as the starting ideal, and with derivative operators $L_{g_1}$ and $L_{g_2}$, where ${g_1}=[x_2~~0]^T$, and $g_2=[0~~x_1^2]^T$. 
\begin{eqnarray} \nonumber
 &&  J_0\coloneqq I_{M_0}=\langle x_1^2x_2 \rangle,  \\ \nonumber
 && L_{g_1}(x_1^2x_2)=2x_1x_2^2 \notin J_0, ~~~~L_{g_2}(x_1^2x_2)=x_1^4\notin J_0, \\ \nonumber
&& \implies J_1\coloneqq J_0 \cup \langle x_1x_2^2, x_1^4 \rangle =  \langle x_1^2x_2,x_1x_2^2, x_1^4 \rangle \\ \nonumber
&& L_{g_1}(x_1x_2^2)=x_2^3\notin J_1,~~~~ L_{g_2}(x_1x_2^2)=2x_1^3x_2\in J_1, \\  \nonumber &&L_{g_1}(x_1^4)=4x_1^3x_2\in J_1,
~~~~L_{g_2}(x_1^4)=0, \\ \nonumber
&&  \implies J_2\coloneqq J_1 \cup   x_2^3   =  \langle x_1^2x_2,x_1x_2^2, x_1^4 , x_2^3\rangle, \\ \nonumber
&& L_{g_1}(x_2^3)=0, ~~~~ L_{g_2}(x_2^3)=3x_1^2x_2^2\in J_2, 
\end{eqnarray}
which shows that  $J_2$ is closed under $L_{g_1}$ and $L_{g_2}$. Therefore $\mathcal{V}(J_2)=(0,0)$ is the only singular point of accessibility.
\end{example}


\subsection{Singular points of strong accessibility}

\begin{theorem} \label{lemma315}
For a generically strongly accessible polynomial system (\ref{sys}), consider the sets $S_{\infty}^*$ and $S_\infty$ as described in Definition \ref{def5}. Then  $S_{\infty}^*=S_{\infty}$, and therefore accessibility from $x$ implies strong accessibility from $x$ and vice versa.
\end{theorem}
\begin{pf} 
From the assumption of analyticity  and generic strong accessibility of the system,  $S_{\infty}^*$ is a closed zero-measure set. 
For each $x_0\in S_{\infty}^*$, let $I(C,x_0)$ be the maximal integral manifold of the accessibility distribution $C$ that passes through $x_0$. By Corollary 3.4 of \cite{sussmann}, the dimension of the strong accessibility distribution $C_0(x)$ is the same for all $x\in I(C,x_0)$, and therefore $I(C,x_0)\subset  S_{\infty}^*$. Moreover $I(C,x_0)$ is a forward-invariant set for the system. 
 Therefore $S_{\infty}^*$ is a zero-measure forward-invariant set for the system, and as a result of Theorem \ref{theorem3}  we have $S_{\infty}^*\subset S_{\infty} $. On the other hand, by definition, we have $C_0(x)\subset C(x)$ for any $x\in D$, which means $S_\infty \subset S_\infty ^*$. From these two inclusion relations we have $S_\infty = S_\infty ^*$. \qed
\end{pf}
As a result of the previous Theorem and Theorem \ref{theorem33}, if the generic rank of $C_0^{n-1}$ is less than $n$, then the system is strongly non-accessible everywhere. Otherwise, the system is generically strongly accessible, and $S_{\infty}^*=S_{\infty}$. Therefore all results of Subsection \ref{3.1} on finding $S_{\infty}$ can be used for finding singular points of strong accessibility. Also since by construction  $C^k$ is generated by the same set of vector fields that generate $C_0^k$ plus $f$,  the strong accessibility index of the system is equal to the accessibility index, or greater by one. 

\subsection{A module-theoretic approach to finding $S_\infty$}
 In what follows, we propose an alternative approach to construct $S_\infty$ that does not need the computation of real radical of ideals.
\begin{lemma} \label{lemma4}
Consider a vector field $Z\in \R^n[x]$, and   a submodule $M$ of the module $\R [x]^n$ over the ring $\R[x]$ that is generated from $\{ X_1,...,X_k \}$. Then  $M$ is invariant under the operator $ad_Z$ iff $ad_Z(X)\in M$ for any $X \in \{X_1,...,X_k \}$. 
\end{lemma}
\begin{pf}
 Since a module contains its generators,  the necessity part is obvious. To prove the sufficiency, note that for every $Y\in M$ there is (at least) one set of $\{ p_1,...,p_k \}$ $\in$ $\R [x]$ such that $Y=\sum_{i=1}^k p_iX_i$. Now, using the properties of Lie bracket stated in (\ref{1a})-(\ref{1c}), the claim can be proved easily. \qed
\end{pf}


\begin{theorem}   \label{theorem4}
For the polynomial system (\ref{sys}),   denote by $C^{\#k}$ a module over $\R [x]$ that is generated by vector fields from $\mathcal{C}^k$. Then there exists an integer $\hat{r}\geq r^*$ such that 
\begin{equation} \label{18}
C^{\#1} \subsetneq C^{\#2} \subsetneq \cdots \subsetneq C^{\#\hat{r}} =C^{\#\hat{r}+1}= \cdots \coloneqq\mathcal{C^{\#}}
\end{equation}
furthermore, $S_{\hat{r}} =S_\infty$. 
\end{theorem}
\begin{pf}
We have $\mathcal{C}^1 \subseteq \mathcal{C}^2 \subseteq \cdots$ by construction, and consequently $C^{\#1} \subseteq C^{\#2} \subseteq \cdots$. This ascending chain of Noetherian modules must stabilize eventually. Assume that $\hat{r}$ is the smallest integer such that $C^{\#\hat{r}} =C^{\#\hat{r}+1}$. By construction, for every $X \in \mathcal{C}^{\hat{r}}$, we have $ad_f(X), ad_{g_i}(X)\in \mathcal{C}^{\hat{r}+1}, ~i=1,...,m$, and therefore $ad_f(X), ad_{g_i}(X)\in C^{\#\hat{r}+1}, ~ i=1,...,m$. Therefore from the assumption $C^{\#\hat{r}} =C^{\#\hat{r}+1}$ and  Lemma \ref{lemma4} we obtain that  $C^{\#\hat{r}}$ is closed under $ad_f, ad_{g_i}$, and because for every $k$ the vector fields of $\mathcal{C}^{k+1}$ are obtained by successive application of operators $ad_f , ad_{g_i}$ on the vector fields of $\mathcal{C}^k$, by a simple induction we obtain $C^{\#\hat{r}} =C^{\#\hat{r}+1}=C^{\#\hat{r}+2}=\cdots$.  Because the columns of each matrix $M_k$ in the proof of Lemma \ref{lemma2} are the generators of $C^{\#{k}}$, from the last equalities  we get that for every $k\geq \hat{r}$, every minor of $M_k$ belongs to the ideal generated by the minors of $M_{\hat{r}}$ of the same dimension, and therefore $S_{\hat{r}} =S_{\hat{r}+1}=\cdots =S_\infty$. Since it was assumed in Lemma \ref{lemma2} that $r^*$ is the smallest integer such that $S_{r^*} =S_\infty$, therefore we have $\hat{r} \geq r^*$. \qed
\end{pf}

\begin{remark}
Theorem \ref{theorem4} suggests that in order to obtain an upper bound on $r^*$, it suffices to look for the first integer $k$ such that two successive submodules generated from $\mathcal{C}^k$ and $\mathcal{C}^{k+1}$, become identical. Identity of two submodules can be checked using the Gr\"{o}bner bases for modules \cite{adams}.
\end{remark}

A similar approach can be taken for determination of singular points of strong accessibility distribution, and therefore we state the following theorem without proof.

\begin{theorem}  \label{theorem4b}
For the polynomial system (\ref{sys}), denote by $C_0^{\#k}$ a module over $\R [x]$ that is generated by vector fields from $\mathcal{C}_0^k$. Then there exists an integer $\hat{l}>l^*$ such that 
\be  \label{19}
C_0^{\#1} \subsetneq C_0^{\#2} \subsetneq \cdots \subsetneq C_0^{\#\hat{l}} =C_0^{\#\hat{l}+1}= \cdots =C_0^{\#}
\ee
furthermore, $S_{\hat{l}}^* =S_{\infty}^*$. 
\end{theorem}

\subsection{Finding singular points with specific rank} \label{3.4}
It may be desirable to find the set of all points $x$ for which $C(x)$ or $C_0(X)$ has dimension less than $l$, for some specific $l<n$. For example, in the case when the system is not generically accessible, one may be interested in finding points at which the rank of  accessibility distribution drops from its generic value. See Examples \ref{ex43} and \ref{ex44}  for the other applications.
For this, we define the set $S_k^{<l }$  as the set of all points $x\in D$ at which  $\textup{dim}\!~C^k(x)<l$, and denote by $S_{\infty}^{<l}$ the set of all points at which  $\textup{dim}\!~C (x)<l$.
Analogously, we define the set $S_{k}^{*<l}$ (respectively $S_{\infty}^{*<l}$) as the set of points at which the distribution $C_0^k$ (respectively $C_0$) has rank less than $l<n$.
 By \cite{sussmann}, $S_{\infty}^{<l}$ (respectively $S_{\infty}^{*<l}$) is the union of all maximal integral manifolds of $C$ (respectively $C_0$) of dimension less than $l$, and therefore the locus of all points $x\in D$ at which 
  $\mathcal{R}_{\Sigma}(x_0)$ (respectively $\mathcal{R}_{\Sigma}(x_0,t)$) has empty interior   in every submanifold of $D$ of dimension greater than or  equal to $l$. \\
 The following theorems show how to obtain $S_{\infty}^{<l}$ and $S_{\infty}^{*<l}$.
 
 \begin{theorem} \label{theorem29}
 The set $S_{\infty}^{<l}$ is an algebraic set,  and  $S_{\infty}^{<l}=S_{\hat{r}}^{<l}$, where $\hat{r}$ is as in (\ref{18}).
 \end{theorem}
 
 \begin{pf}
Assume that $\{b_{k,1},...,b_{k,p}\}$ are all $l\times l$ minors of $M_k$, and define $I_{M_k}^{<l}\coloneqq \langle b_{k,1},...,b_{k,p} \rangle $. Then we have $S_k^{<l}=\mathcal{V} (I_{M_k}^{<l})$ by construction. Similarly to (\ref{4}), we have the following descending chain of algebraic sets
 \begin{equation}  \label{20}
S_1^{\tiny <l} \supseteq S_2^{<l} \supseteq \cdots  \supseteq S_{r_l^*}^{<l} =S_{r_l^*+1}^{<l}= \cdots = S_{\infty}^{<l}
\end{equation}
that eventually stabilizes at the algebraic set  $S_{\infty}^{<l}$.
Also,  as it was shown in the proof of Theorem \ref{theorem4}, for  $k>\hat{r}$ the ascending chain of modules $C^{\#k}$ stabilizes, and every $l\times l$ minor of  $M_k$ belongs to the ideal $I_{M_{\hat{r}}}^{<l}$ and therefore $S_{\infty}^{<l}=S_{\hat{r}}^{<l}=\mathcal{V}(I_{M_{\hat{r}}}^{<l})$.
  \qed
 \end{pf}


  \begin{theorem}
  For a generically strongly accessible polynomial system (\ref{sys}), consider the sets $S_{\infty}^{*<l}$ and $S_{\infty}^{<l}$ as described in Definition \ref{def5}. Then  $S_{\infty}^{*<l}=S_{\infty}^{<l}$. Furthermore,  $S_{\infty}^{*<l}=S_{\hat{l}}^{*<l}$, where $\hat{l}$ is as in (\ref{19}).
\end{theorem}
\begin{pf}
The proof of $S_{\infty}^{*<l}=S_{\infty}^{<l}$ is similar to the proof of $S_{\infty}^{*}=S_{\infty}$ in Theorem \ref{lemma315}.
The proof of the last part of the theorem is similar to the proof of the last part of Theorem \ref{theorem29}, except that the modules $C_{0}^{\#k}$, the sets $S_k^{*<l}$ and the $l \times l$ minors of $M_k^*$ must be considered.  \qed
\end{pf}

\section{Non-polynomial systems} \label{sec5}

For an input-affine system with analytic or smooth vector fields $\{ f,g_1,...,g_m\}$, it is still possible to define matrices $M_k$ and correspondingly ideals $I_{M_k}$ and the sets $S_k$, albeit in a non-Noetherian ring. So, as the Hilbert Basis Theorem doesn't hold in a non-Noetherian ring, one may find examples to show that in smooth or analytic systems,  the descending chain of sets $S_k$ in (\ref{3}) may  never stabilize, and therefore the accessibility index of the system be $\infty$.
\begin{example}

Let $\alpha (x_1)\coloneqq \prod_{n=1}^{\infty}(1-\frac{ x_1}{a_n})^n $ where the constants $a_n$ being chosen in such a way that the product is convergent for any $x$. Let the system $\Sigma$ be given by $\dot{x}_1=u,~ \dot{x}_2=\alpha (x_1)$. Then $S_k= \{x\in \R^2 | x_1=a_{n+1}, ~n\geq k\}$. Thus the sequence $S_k$ does not stabilize, and $r^*=\infty$.
\end{example} 

Fortunately, on compact semianalytic sets, the descending chain property holds for any chain of analytic sets, and therefore many of the results can be extended to the case of analytic systems. But, to keep things simple, and also to take advantage of the  computational power of computer algebra tools, here we only consider the wide class of analytic nonlinear control systems that can be simplified to polynomial form, using the immersion technique \cite{zbig, ohtsuka, carravetta}. System immersion  is usually performed by defining some functions of $x$ as new state variables, and  may cause an increase in the dimension of the system. 
 The system (\ref{sys}) is said to be (invariantly) immersible \cite{ohtsuka} into a polynomial system, if there exist an analytic immersion mapping $z\coloneqq T(x) : \R^n \rightarrow \R^{n^*}$,  and an  $n^*$-dimensional polynomial system 
 $  \label{immersedsys}
  \hat{\Sigma}_z: \dot{z}=\hat{f}(z) + \Sigma_{i=1}^m u_i \hat{g}_i(z),
   $
  where $\hat{f},\hat{g_i}$ are push-forwards of $f,g_i$, respectively, by $T(x)$, (i.e.  $L_f T(x)=\hat{f}(T(x))$,  $L_{g_i} T(x)=\hat{g_i}(T(x))$).
Recall that a mapping $ T(x) : \R^n \rightarrow \R^{n^*}$ is called a (local) immersion if $\textup{dim}\!~dT(x)=n$ for every $x$, so $n^* \geq n$.\\ 
 For the analytic system (\ref{sys}), denote by $\mathcal{O}$ the smallest subspace (over $\R$) of $\mathcal{A}$ that contains $\{x_1,\dots ,x_n\}$ and is invariant under $L_f,L_{g_1},\dots, L_{g_m}$.
  Using the results of \cite{ohtsuka}, a sufficient condition for immersibility into a polynomial system is that $\mathcal{O}$ be a subset of a finitely generated field over $\R$. 
For example, for the system $\Sigma : \dot{x}_1=u\sin(x_2),~~ \dot{x}_2=x_1$,  the set $\mathcal{O}$ is a subset of the field that is generated by $\{x_1,x_2,\sin(x_2),\cos(x_2)\}$ over $\R$, and therefore  the immersion mapping  $z= T(x)\coloneqq (x_1,x_2,\sin(x_2),\cos(x_2))$ transforms the system into the four-dimensional polynomial system $\hat{\Sigma}: \dot{z}_1=uz_3, \dot{z}_2=z_1, \dot{z}_3=uz_3z_4, \dot{z}_4=-z_3^2$.
  
  Denote by $\hat{S}_{\infty}^{< n}$ (respectively $\hat{S}_{\infty}^{*< n}$) the set of points of $\R^{n^*}$ at which the accessibility (respectively strong accessibility) distribution of $\hat{\Sigma}$ has rank less than $n$. Then the following theorem relates the singular points of accessibility of the original system (\ref{sys}) to $\hat{S}_{\infty}^{< n}$ and $\hat{S}_{\infty}^{*< n}$.
  \begin{theorem} \label{theorem42}
  Assume that the system (\ref{sys}) is immersible into a polynomial system $\hat{\Sigma}$ by an immersion mapping $T : \R^n \rightarrow \R^{n^*}$.   Then for every $x\in \R^n$, we have $x\in S_{\infty}$ (respectively $x\in S_{\infty}^*$) iff $T(x)\in \hat{S}_{\infty}^{< n}$ (respectively $T(x)\in \hat{S}_{\infty}^{*< n}$) and therefore accessibility index (respectively strong accessibility index) of the system (\ref{sys}) is finite.
  \end{theorem}

 \begin{pf}
 Without loss of generality,  we assume that $T_i(x)=x_i$, for $ 1\leq i \leq n$. Then   $T : \R^n \rightarrow \R^{n^*}$ is an injective immersion, and hence the mapping $T$ is an embedding of the manifold $\R^n$ into the image of $T$, everywhere diffeomorphic,  and the image of $T$ is an $n$-dimensional invariant manifold for the system $\hat{\Sigma}$. Denote by $\hat{C^k}$ (respectively $\hat{C_0^k}$) the accessibility distribution of order $k$ (respectively strong accessibility distribution of order $k$) of the system $\hat{\Sigma}$. 
   Since $\hat{f},\hat{g}_1,\dots,\hat{g}_m$ are push-forwards of $f,g_1,\dots,g_m$ by the diffeomorphism $T$, and  push-forward of vector fields by a diffeomorphism commutes with Lie bracketing, for any $x_0\in \R^n$ we have $\textup{dim}(\hat{C}^k(T(x_0)))=\textup{dim}({C}^k(x_0))$ (respectively $\textup{dim}(\hat{C}_0^k(T(x_0)))=\textup{dim}({C}_0^k(x_0))$), and therefore the claim follows easily. \qed
   \end{pf}
      The following  examples demonstrates the application of Theorem \ref{theorem42} in non-polynomial systems.
  \begin{example} \label{ex43}
  We test the equations of a unicycle, for possible singular points of accessibility distribution. The dynamics of the unicycle system is described by $\Sigma: \dot{x}_1=u_1\cos(x_3), ~\dot{x}_2=u_1\sin (x_3), ~\dot{x}_3=u_2$.  
 The transformation $z = T(x):\R ^3 \rightarrow \R ^5$ defined by $T(x) \coloneqq (x_1 , x_2,x_3,$ $\sin (x_3),$ $\cos (x_3))^T $ immerses the system $\Sigma$ into the  polynomial system
$\hat{\Sigma}: \dot{z}=u_1 \hat{g_1}(z) + u_2  \hat{g_2}(z)$,
  with $\hat{g_1}(z)\coloneqq [z_5,z_4,0,0,0]^T$ and $ \hat{g_2}(z)\coloneqq [0, 0, 1,z_5 , -z_4 ]^T$.\\   
  The set $S_\infty$ of the system $\Sigma$, corresponds to the intersection of the set $\hat{S}_{\infty}^{<3}$ of the system $\hat{\Sigma}$ with   $\textup{im}~\!T$. We use the results of Section \ref{3.4} to obtain $\hat{S}_{\infty}^{<3}$. Computing $[\hat{g_1},\hat{g_2}]=[ z_4 ~~-z_5~~ 0 ~~ 0 ~~0 ]^T$,  $ [\hat{g_1},[\hat{g_1},\hat{g_2}]]=[ 0 ~~0~~0 ~~ 0 ~~ 0]^T$, $[\hat{g_2},[\hat{g_1},\hat{g_2}]]=[ z_5~~z_4~~ 0 ~~0 ~~ 0 ]^T=\hat{g_1}$ 
  shows that $\hat{C}_1^\#$, which is the module over $\R [z]$ generated by $\hat{g}_1,~ \hat{g}_2, ~ [\hat{g}_1,\hat{g}_2]$, is invariant under $ad_{\hat{g}_1}$ and $ad_{\hat{g}_2}$, and therefore $\hat{S}_{\infty}^{<3}=\mathcal{V}(I_{M_1}^{<3})$. Computing $M_1=[ \hat{g}_1~~ \hat{g}_2~~ [\hat{g}_1,\hat{g}_2] ]$, and $I_{M_1}^{<3}$ as the ideal generated by all $3 \times 3$ minors of $M_1$, we have $I_{M_1}^{<3}=\langle z_4^2 +z_5^2\rangle$. So $\hat{S}_{\infty}^{<3}=\{ z\in \R^5~ |~ z_4^2 +z_5^2=0 \}$, and the intersection of $\hat{S}_{\infty}^{<3}$  by the set $\textup{im}~\!T=\{z\in \R^5 ~|~ z_4=\sin(z_3), ~ z_5=\cos(z_3) \}$ is empty set. Therefore $\Sigma$ is accessible everywhere. 
  \end{example}  
  \begin{example} \label{ex44}
  Consider a vertically driven pendulum system $\Sigma$, where the base is constrained to move only vertically on a bar, the base and the pendulum both of unit mass, and unit length. The vertical input force $u$ acts on the base of pendulum, with positive $u$ for upward forces. This system has four (independent) state variable $y, \dot{y}, \theta, \dot{\theta} $ , where $y$ is the place of the base on the bar, and $\theta$ is the angle between the pendulum and the bar. The equations of motion are
\begin{eqnarray}
&&  \hskip -0.7cm\nonumber  \ddot{y}=-10+u-\cos{(\theta)}\dot{\theta}^2 \frac{1}{2-\sin^2{(\theta)}} +u \frac{1}{2-\sin^2{(\theta)}}\\
&&  \hskip -0.7cm \nonumber   \ddot{\theta}=\sin{(\theta)} \cos{(\theta)}\dot{\theta}^2 \frac{1}{2-\sin^2{(\theta)}} +u sin{(\theta)} \frac{1}{2-\sin^2{(\theta)}}.
    \end{eqnarray}  
    The transformation $z = T(y, \dot{y}, \theta, \dot{\theta}):\R ^4 \rightarrow \R ^7$ defined by $T(y, \dot{y}, \theta, \dot{\theta}) \coloneqq (y , \dot{y},\theta, \dot{\theta}, \sin(\theta), \cos(\theta), \frac{1}{2-\sin^2{(\theta)}})^T $ immerses the system $\Sigma$ into the  polynomial system
 $\hat{\Sigma}: \dot{z}=\hat{f}(z) + u  \hat{g}(z)$
  with $\hat{f}(z)\coloneqq [z_2,~ z_4^2z_6z_7-10,~ z_4,$ $~ z_4^2z_5z_6z_7,$ $ ~z_4z_6,$ $~ -z_4z_5,~ 2z_4z_5z_6z_7^2]^T$ and $ \hat{g}(z)\coloneqq [0,~ z_7,~ 0,~z_5z_7 ,$ $~ 0,~0,~0 ]^T$.  
  First we obtain $S_{\infty}^{<4}$ for $\hat{\Sigma}$  by use of the results of Section \ref{3.4}. Computing the vector fields of distributions $\hat{C}_k$, and using Gr\"{o}bner bases for modules, we check whether the Module  $C_k^\#$ is invariant under $ad_{\hat{f}}$ and $ ad_{\hat{g}}$ or not. After a few computations, it turns out that the ascending chain of modules finally stabilizes at $C_6^\#$. Computing $S_{\infty}^{<4}=\mathcal{V}(I_{M_6}^{<4})=\{ z\in \R^7 ~|~ z_4z_6z_7\!=\!0~ \&~ z_5z_7\!=\!0\}$, the set of accessibility singular points of the system $\Sigma$ correspond to $S_{\infty}^{<4}\cap \textup{im}~\!T= \{ z\in \R^7 | z_4\!=\!0~ \&~ z_5=0\}$, which in the coordinates of the original system $\Sigma$ reads as $\dot{\theta}=0$ and $\sin{\theta}=0$.
  \end{example} 

\section{Conclusion}

The paper addresses the problem of finite determination of accessibility/strong accessibility for two large subclasses of nonlinear systems, namely polynomial systems and analytical systems that are immersible into the polynomial systems. 
 It is shown that the set of accessibility singular points is the maximal zero-measure invariant set of the system. Thanks to  the descending chain property and invariance of this set, several theorems and algorithms are stated to obtain the entire set of singular points, as well as the minimum number of lie brackets in the accessibility rank   test that is necessary  for deciding accessibility from any point, called accessibility index in this paper.  Alternative algorithms are proposed that compute upper bounds on accessibility index that are easier to find. The solved real-life examples shows the applicability of the results using computer algebra tools, an improvements over the previously obtained bounds. 

%
%
%
%

%

\section*{Acknowledgment}
The work of Zbigniew Bartosiewicz has been supported by grant No.S/WI/1/2016 of Bialystok University of Technology, financed by Polish Ministry of Science and Higher Education.

\appendix
\section{Appendix}
We recall some basic facts from  real algebraic geometry.  A subset   $I\subset \R[x]$ is an ideal of $\R[x]$ if for any $a,b\in I$ and $c\in \R[x]$ we have  $a+b\in I$, and $ca\in I$. An ideal is said to be proper if it does not contain 1. For a given set of polynomials $p_1,...,p_r \in \R [x]$, the ideal generated by $p_1,...,p_r$ is defined as
\[ \langle p_1,...,p_r \rangle \coloneqq      \Big\{ \sum_{i=1}^{r} a_i p_i : a_1,...,a_r \in \R [x]  \Big\}. \]
Since $\R [x] $ is a Noetherian ring, every ideal $I \in \R [x] $ is finitely generated \cite{bochnak}.

For an ideal $I$ of $\R[x]$, its radical, denoted by $\sqrt{I}$, is the set of all $p\in \R [x]$ such that $p^k \in I$ for some $k \in \N$. The real radical of $I$, denoted by $\sqrt[\R]{I}$, is the set of all $p \in \R [x]$ for which there exist $q_1,...,q_k \in \R [x]$ and $m,k \in \N$, such that $p^{2m} +\sum_{i=1}^{k} q^2 \in I$.

\begin{proposition} \label{prop1}
\cite{bochnak} If $I$ and $J$ are ideals of $\R [x]$, then the following holds
\begin{enumerate}
\item [(i)] The real radical of $I$ is an ideal of $\R [x]$,
\item [(ii)] $I \subseteq \sqrt{I} \subseteq \sqrt[\R]{I},$
\item [(iii)] if $I \subseteq J$ then $ \sqrt[\R]{I} \subseteq \sqrt[\R]{J}$.   \color{red}
\end{enumerate}
\end{proposition}

The  algebraic set of an ideal $I \subset \R [x]$ is defined as $\mathcal{V}(I) \!\coloneqq  \!\{ x\in \R^n : p(x)\!=\!0 \textup{ for all } p\in I \}$. In other words, $\mathcal{V}(I)$ is the set of common zeroes of all polynomials in $I$.

For a subset $A \subset \R^n$, its zero-ideal, denoted by $\mathcal{I}(A)$, is defined as  $ \mathcal{I}(A) =\{ p\in \R [x] : p(x)=0 \textup{ for all } x\in A \}$.

\begin{proposition}  \label{prop2}
\cite{bochnak} Let $I$ be an ideal of $\R [x]$. Then $\mathcal{I}( \mathcal{V} (I))$ $=$ $\sqrt[\R]{I}$.
\end{proposition}
\begin{proposition} \label{prop3}
\cite{bochnak} Let $I$ and $J$ be  ideals of $\R [x]$,  and $A$ and $B$  subsets of $\R^n$. Then
\begin{enumerate}
\item [(i)] $ A \subseteq \mathcal{V}(\mathcal{I}(A))$
\item [(ii)] $I \subseteq \mathcal{I}(\mathcal{V}(I))$
\item [(iii)] $\mathcal{V}(\mathcal{I}(\mathcal{V}(I)))= \mathcal{V}(I)$
\item [(iv)] $\mathcal{I}(\mathcal{V}(\mathcal{I}(A)))=\mathcal{I}(A)$ 
\item [(v)] if $I \subseteq J$ then $\mathcal{V}(J) \subseteq  \mathcal{V}(I)$.
\item [(vi)] if $A\subseteq B$ then $\mathcal{I}(B) \subseteq \mathcal{I}(A)$.
\end{enumerate}
\end{proposition}


\begin{thebibliography}{21}
\bibitem[{Sussmann and Jurdjevic(1972)}]{sussmann}
\bibinfo{author}{H.~J. Sussmann}, \bibinfo{author}{V.~Jurdjevic},
\newblock \bibinfo{title}{Controllability of nonlinear systems},
\newblock \bibinfo{journal}{Journal of Differential Equations}
  \bibinfo{volume}{12} (\bibinfo{year}{1972}) \bibinfo{pages}{95--116}.
\bibitem[{Kawski(2006)}]{kawski}
\bibinfo{author}{M.~Kawski},
\newblock \bibinfo{title}{On the problem whether controllability is finitely
  determined},
\newblock in: \bibinfo{booktitle}{Proceedings of MTNS},
  volume~\bibinfo{volume}{6}, \bibinfo{year}{2006}.
\bibitem[{Gabrielov and Khovanskii(1998)}]{gabriel}
\bibinfo{author}{A.~Gabrielov}, \bibinfo{author}{A.~Khovanskii},
  \bibinfo{title}{Multiplicity of a {Noetherian} intersection}, volume
  \bibinfo{volume}{186}, \bibinfo{publisher}{Amer. Math. Soc., Providence, RI},
  \bibinfo{year}{1998}, pp. \bibinfo{pages}{119--130}.
\bibitem[{Gabrielov(1995)}]{gabrielov}
\bibinfo{author}{A.~Gabrielov},
\newblock \bibinfo{title}{Multiplicities of zeroes of polynomials on
  trajectories of polynomial vector fields and bounds on degree of
  nonholonomy},
\newblock volume~\bibinfo{volume}{2}, \bibinfo{publisher}{Amer. Math. Soc.,
  Providence, RI}, \bibinfo{year}{1995}, pp. \bibinfo{pages}{437--451}.
\bibitem[{Risler(1996)}]{risler}
\bibinfo{author}{J.-J. Risler},
\newblock \bibinfo{title}{A bound for the degree of nonholonomy in the plane},
\newblock \bibinfo{journal}{Theoretical Computer Science} \bibinfo{volume}{157}
  (\bibinfo{year}{1996}) \bibinfo{pages}{129--136}.
\bibitem[{Gabrielov(1999)}]{gabrielov2}
\bibinfo{author}{A.~Gabrielov},
\newblock \bibinfo{title}{Multiplicity of a zero of an analytic function on a
  trajectory of a vector field},
\newblock in: \bibinfo{booktitle}{The Arnoldfest (Toronto, ON, 1997)},
  \bibinfo{publisher}{American Mathematical Society},
  \bibinfo{address}{Providence, RI}, \bibinfo{year}{1999}, pp.
  \bibinfo{pages}{191--200}.
\bibitem[{Binyamini(2016)}]{binyamini}
\bibinfo{author}{G.~Binyamini},
\newblock \bibinfo{title}{Multiplicity estimates: a {M}orse-theoretic
  approach},
\newblock \bibinfo{journal}{Duke Mathematical Journal} \bibinfo{volume}{165}
  (\bibinfo{year}{2016}) \bibinfo{pages}{95--128}.
\bibitem[{M{\"u}ller and Donelan(2018)}]{muller}
\bibinfo{author}{A.~M{\"u}ller}, \bibinfo{author}{P.~Donelan},
\newblock \bibinfo{title}{Towards a unified notion of kinematic singularities
  for robot arms and non-holonomic platforms},
\newblock in: \bibinfo{booktitle}{Advances in Robot Kinematics 2016},
  \bibinfo{publisher}{Springer}, \bibinfo{year}{2018}, pp.
  \bibinfo{pages}{393--401}.
\bibitem[{Tcho{\'n}(2000)}]{tchon}
\bibinfo{author}{K.~Tcho{\'n}},
\newblock \bibinfo{title}{On kinematic singularities of nonholonomic robotic
  systems},
\newblock in: \bibinfo{booktitle}{Romansy 13}, \bibinfo{publisher}{Springer},
  \bibinfo{year}{2000}, pp. \bibinfo{pages}{75--84}.
\bibitem[{Kaminski et~al.(2018)Kaminski, L{\'e}vine, and Ollivier}]{kaminski}
\bibinfo{author}{Y.~J. Kaminski}, \bibinfo{author}{J.~L{\'e}vine},
  \bibinfo{author}{F.~Ollivier},
\newblock \bibinfo{title}{Intrinsic and apparent singularities in
  differentially flat systems, and application to global motion planning},
\newblock \bibinfo{journal}{Systems \& Control Letters} \bibinfo{volume}{113}
  (\bibinfo{year}{2018}) \bibinfo{pages}{117--124}.
\bibitem[{Sarafrazi et~al.(2019)Sarafrazi, Pawluszewicz, Bartosiewicz, and
  Kotta}]{amin}
\bibinfo{author}{M.~A. Sarafrazi}, \bibinfo{author}{E.~Pawluszewicz},
  \bibinfo{author}{Z.~Bartosiewicz}, \bibinfo{author}{{\"U}.~Kotta},
\newblock \bibinfo{title}{On the finiteness of accessibility test for nonlinear
  discrete-time systems},
\newblock \bibinfo{journal}{arXiv preprint, arXiv:1905.10154}
  (\bibinfo{year}{2019}).
\bibitem[{Zariski and Samuel(2013)}]{zariski}
\bibinfo{author}{O.~Zariski}, \bibinfo{author}{P.~Samuel},
  \bibinfo{title}{Commutative Algebra}, volume~\bibinfo{volume}{2},
  \bibinfo{publisher}{Springer Science \& Business Media},
  \bibinfo{year}{2013}.
\bibitem[{Nijmeijer and van~der Schaft(1990)}]{nijmeijer}
\bibinfo{author}{H.~Nijmeijer}, \bibinfo{author}{A.~van~der Schaft},
  \bibinfo{title}{Nonlinear Dynamical Control Systems},
  \bibinfo{publisher}{Springer New York}, \bibinfo{year}{1990}.
\bibitem[{Sussmann(1973)}]{sussmann2}
\bibinfo{author}{H.~J. Sussmann},
\newblock \bibinfo{title}{Orbits of families of vector fields and integrability
  of distributions},
\newblock \bibinfo{journal}{Transactions of the American Mathematical Society}
  \bibinfo{volume}{180} (\bibinfo{year}{1973}) \bibinfo{pages}{171--188}.
\bibitem[{Nagano(1966)}]{nagano}
\bibinfo{author}{T.~Nagano},
\newblock \bibinfo{title}{Linear differential systems with singularities and an
  application to transitive lie algebras},
\newblock \bibinfo{journal}{Journal of the Mathematical Society of Japan}
  \bibinfo{volume}{18} (\bibinfo{year}{1966}) \bibinfo{pages}{398--404}.
\bibitem[{Cox et~al.(2007)Cox, Little, and O'Shea}]{cox}
\bibinfo{author}{D.~Cox}, \bibinfo{author}{J.~Little},
  \bibinfo{author}{D.~O'Shea}, \bibinfo{title}{Ideals, Varieties, and
  Algorithms}, volume~\bibinfo{volume}{3}, \bibinfo{publisher}{Springer},
  \bibinfo{year}{2007}.
\bibitem[{Adams and Loustaunau(1994)}]{adams}
\bibinfo{author}{W.~W. Adams}, \bibinfo{author}{P.~Loustaunau},
  \bibinfo{title}{An Introduction to Gr\"{o}bner Bases},
  \bibinfo{publisher}{American Mathematical Soc.}, \bibinfo{year}{1994}.
\bibitem[{Bartosiewicz(1986)}]{zbig}
\bibinfo{author}{Z.~Bartosiewicz},
\newblock \bibinfo{title}{Realizations of polynomial systems},
\newblock in: \bibinfo{editor}{M.~Fliess}, \bibinfo{editor}{M.~Hazewinkel}
  (Eds.), \bibinfo{booktitle}{Algebraic and geometric methods in nonlinear
  control theory}, \bibinfo{publisher}{Springer Netherlands},
  \bibinfo{address}{Dordrecht}, \bibinfo{year}{1986}, pp.
  \bibinfo{pages}{45--54}.
\bibitem[{Ohtsuka(2005)}]{ohtsuka}
\bibinfo{author}{T.~Ohtsuka},
\newblock \bibinfo{title}{Model structure simplification of nonlinear systems
  via immersion},
\newblock \bibinfo{journal}{IEEE Transactions on Automatic Control}
  \bibinfo{volume}{50} (\bibinfo{year}{2005}) \bibinfo{pages}{607--618}.
\bibitem[{Carravetta(2015)}]{carravetta}
\bibinfo{author}{F.~Carravetta},
\newblock \bibinfo{title}{Global exact quadratization of continuous-time
  nonlinear control systems},
\newblock \bibinfo{journal}{SIAM Journal on Control and Optimization}
  \bibinfo{volume}{53} (\bibinfo{year}{2015}) \bibinfo{pages}{235--261}.
\bibitem[{Bochnak et~al.(2013)Bochnak, Coste, and Roy}]{bochnak}
\bibinfo{author}{J.~Bochnak}, \bibinfo{author}{M.~Coste},
  \bibinfo{author}{M.-F. Roy}, \bibinfo{title}{Real Algebraic Geometry},
  volume~\bibinfo{volume}{36}, \bibinfo{publisher}{Springer Science \& Business
  Media}, \bibinfo{year}{2013}.

\end{thebibliography}
\end{document}